\documentclass[11pt]{article}
\usepackage{amssymb,amsmath,amsfonts}
\usepackage{amssymb, amsmath, url,graphicx}

\begin{document}

\author{Mircea Crasmareanu}
\title{A new Tzitzeica hypersurface}
\date{May 12, 2010}
\maketitle

\begin{abstract}
A new hypersurface of Tzitzeica type is obtained in all three forms: parametric, implicit and explicit. Its two-dimensional version, although well-known from a theoretical point of view, is plotted with Matlab.
\end{abstract}



\medskip

%\section*{Introduction}

The first centroaffine invariant was introduced in 1907 by G. Tzitzeica  in the classical
theory of surfaces. Namely, if $M^2$ is such a surfaces embedded in $\mathbb{R}^3$ the
Tzitzeica invariant is $Tzitzeica(M^2)=\frac{K}{d^4}$ where $K$ is the Gaussian curvature and
$d$ is the distance from the origin of $\mathbb{R}^3$ to the tangent plane in an arbitrary point of
$M^2$.

Since then, the class of surfaces with a constant $Tzitzeica$ was the subject of several fruitful research, see for example \cite{b:v} and the surveys \cite{v:c} and \cite{p:g}. For example,
are called {\it Tzitzeica surfaces} by Romanian geometers, {\it affine spheres} by Blaschke and
{\it projective spheres} by Wilczynski. One interesting direction of study is the generalization to
higher dimension and then a hypersurface $M^n$ of $\mathbb{R}^{n+1}$ was called {\it Tzitzeica hypersurface} if
$Tzitzeica(M^n)=\frac{K}{d^{n+2}}$ is a real constant.

The spheres and the quadrics are the simplest
examples of Tzitzeica surfaces. Tzitzeica himself obtains the surface $xyz=1$ which was generalized by Calabi to
$x^1...x^{n+1}=1$ in \cite{e:c}.

The aim of this note is to derive the $n$-dimensional version of another well-known Tzitzeica surface:
$$
z(x^2+y^2)=1 \eqno(0)
$$
The importance of this surface is pointed out by the first theorem of \cite{d:v}.

Inspired by \cite{b:p} we search for $M^n$ being a hypersurface of rotation in the parametrical form:
$$
M^n: x^1=\frac{2u^nu^i}{\Delta },..., x^{n-1}=\frac{2u^nu^{n-1}}{\Delta }, x^n=\frac{u^n(\Delta -2)}{\Delta }, x^{n+1}=f(u^n)
$$
whit $f:(0, +\infty )\rightarrow \mathbb{R}$ and $\Delta =(u^1)^2+...+(u^n)^2+1$.

For easy computations we derive an implicit equation of $M^n$. More precisely, from the above equations we have: $
(u^n)^2=(x^1)^2+...+(x^n)^2$ and then $M^n: x^{n+1}=f(\sqrt{(x^1)^2+...+(x^n)^2})$ i.e.:
$$
M^n: F(x^1,...,x^{n+1})=f(\sqrt{(x^1)^2+...+(x^n)^2})-x^{n+1}=0. \eqno(1)
$$

Recall that choosing the normal $N=-\frac{\nabla F}{\|\nabla F\|}$ the Gaussian curvature of $M^n$ is:
$$
K=-\frac{\left|
\begin{array}{ll}
F_{ij} & F_{i} \\
F_{j} & 0
\end{array}
\right| }{\parallel \nabla F\parallel ^{n+2}}.
$$

For $p=(x^i)\in M^n$ the tangent hyperplane $T_pM$ has the equation \\ $F_i(X^i-x^i)=0$ and then:
$$
d=\frac{|F_ix^i|}{\|\nabla F\|}.
$$

From the last two relations the Tzitzeica condition reads:
$$
\left|
\begin{array}{ll}
F_{ij} & F_{i} \\
F_{j} & 0
\end{array}
\right| =-Tzitzeica(M^n)|F_ix^i|^{n+2}. \eqno(2)
$$

With $F$ of $(1)$ we obtain that the LHS of $(2)$ is $\frac{f''(f')^{n-1}}{\Delta ^{\frac{n-1}{2}}}$
while the RHS of $(2)$ is $-Tzitzeica(M^n)|f'\sqrt{\Delta }-f|^{n+2}$. Therefore, denoting $\sqrt{\Delta }=t$ and considering $f=f(t)$
we have:
$$
f''(f')^{n-1}=-Tzitzeica(M^n)t^{n-1}|tf'-f|^{n+2}
$$
for which we search $f=t^a$. The comparison of degrees of $t$ gives:
$$
a-2+(a-1)(n-1)=n-1+a(n+2)
$$
with solution $a=-n$.

Our result is:

\medskip

{\bf Theorem} {\it The hypersurface of $\mathbb{R}^{n+1}$:
$$
M^n: (x^{n+1})^2[(x^1)^2+...+(x^n)^2]^n=1 \eqno(3)
$$
is a Tzitzeica one with}:
$$
Tzitzeica(M^n)=\frac{(-n)^n}{(n+1)^{n+1}}. \eqno(4)
$$

\medskip

{\bf Remarks} \\
1) If in computing $d$ we choose the opposite normal $-N$ the we obtain $-Tzitzeica(M^n)$.
For example in \cite[p. 137]{u:b} is $M^2$ with $Tzitzeica(M^2)=-\frac{4}{27}$. \\
2) The main theorem of \cite{d:v} gives the 3-dimensional {\it affine locally strongly convex}
hypersurface: $(y^2-z^2-w^2)^3x^2=1$ as generalization of $(0)$. The authors continue their
study in \cite{d:vy} where Theorem 1 states the 5-dimensional variant $(y^2-z^2-w^2-t^2)^2x=1$.

\medskip

With Matlab command: \\
$>>$ ezsurf('$1/(u^2+v^2)$',[-5,5],[-5,5])

\noindent the following picture of $M^2$ is obtained:

\begin{figure}[ht]\label{Fig1}
\centerline{\includegraphics[width=5in]{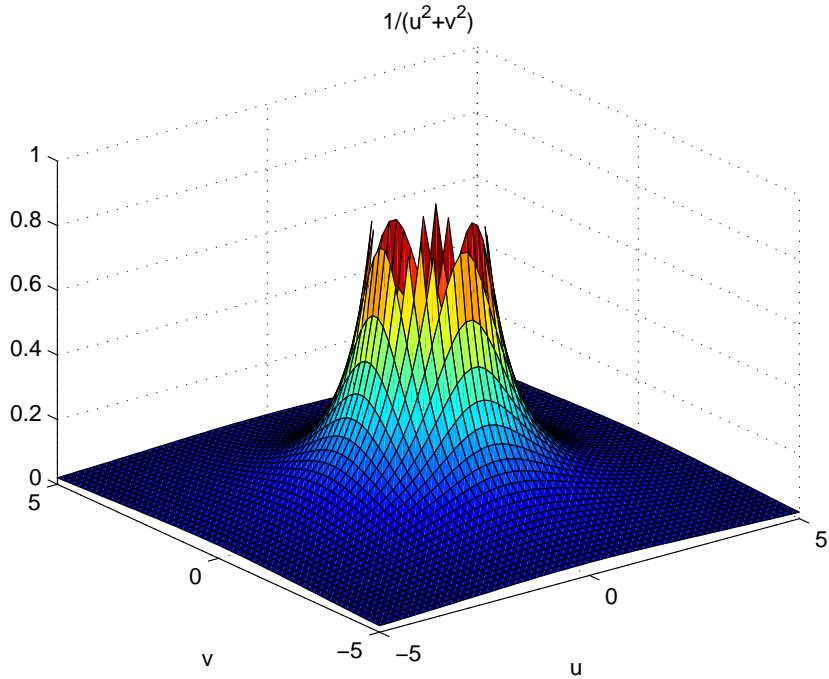}}
\caption{$z(x^2+y^2)=1$.}
\end{figure}

{\bf Acknowledgement} Special thanks are offered to ??? and
% for several useful remarks.

\vspace{.2cm}

\noindent  Faculty of Mathematics \newline
University "Al. I. Cuza" \newline
 Ia\c{s}i, 700506 \newline
 Rom\^{a}nia  \newline
mcrasm@uaic.ro

\smallskip

\noindent http://www.math.uaic.ro/$\sim$mcrasm

\end{document}